\documentclass{ifacconf}

\usepackage{graphicx}      
\usepackage{natbib}        
\usepackage{amsmath,amssymb}

\begin{document}
\begin{frontmatter}

  \title{ Parallel iterative methods for variational integration applied to navigation problems\thanksref{footnoteinfo}} 

\thanks[footnoteinfo]{\copyright 2021 the authors. This work has been accepted to IFAC for publication under a Creative Commons Licence CC-BY-NC-ND. D. Mart{\'\i}n de Diego acknowledges financial support from the Spanish Ministry of Science and Innovation, under grant PID2019-106715GB-C21 and  the ``Severo Ochoa Programme for Centres of Excellence'' in R\&D  and I-Link Project (Ref: linkA20079) from CSIC (CEX2019-000904-S). S.\ Ferraro acknowledges financial support from PICT 2019-00196, FONCyT, Argentina, and PGI 2018, UNS.}

\author[First]{Sebasti\'an J. Ferraro} 
\author[Second]{David Mart\'in de Diego} 
\author[Third]{Rodrigo T. Sato Mart\'in de Almagro}

\address[First]{Instituto de Matem\'atica (INMABB) -- Departamento de Matem\'atica, Universidad Nacional del Sur (UNS) -- CONICET, Bah\'ia Blanca, Argentina (e-mail: sferraro@uns.edu.ar)}
\address[Second]{Instituto de Ciencias Matem\'aticas,  ICMAT
c/ Nicol\'as Cabrera, n$^\textrm{o}$~13-15, Campus Cantoblanco, UAM
28049 Madrid, Spain (e-mail: david.martin@icmat.es)}
\address[Third]{Institute of Applied Dynamics, Friedrich-Alexander-Universit\"at Erlangen-N\"urnberg, Germany (e-mail: rodrigo.t.sato@fau.de)}

\begin{abstract}                
Discrete variational methods have shown an excellent performance in numerical simulations of different mechanical systems. In this paper, we introduce an iterative method for discrete variational methods appropriate for boundary value problems. More concretely, we explore a parallelization strategy that leverages the power of multicore CPUs and GPUs (graphics cards). We study this parallel method for first-order and second-order Lagrangians and we illustrate its excellent behavior in some interesting applications, namely Zermelo's navigation problem, a fuel-optimal navigation problem, and an interpolation problem. 
\end{abstract}


\end{frontmatter}

\section{Introduction}

Variational integrators are numerical algorithms obtained from a discrete variational principle \citep{Marsden_West:Discrete_mechanics_and_variational_integrators}. These methods exhibit excellent structure-pre\-serv\-ing properties as a consequence of their variational derivation. Furthermore, these variational methods preserve qualitative properties of a system such as symmetries, constants of motion, or the manifold structure of the configuration space. They also allow for the preservation of geometric structures such as symplectic or contact structures, or Poisson brackets. Moreover, there exist extensions to other cases of interest such as systems with external forces, holonomic and nonholonomic constraints, optimal control theory and classical field theories.

\subsection{Variational discrete equations}
For a discrete mechanical system on a configuration manifold $Q$, with a discrete Lagrangian $L_d\colon Q\times Q \to \mathbb{R}$, a sequence $\{q_k\}_{k=0}^{N}$ in $Q$ is a trajectory if and only if it satisfies the discrete Euler--Lagrange (\textbf{DEL}) equations
\begin{equation}\label{eq:DEL}
	D_2L_d(q_{k-1}, q_k)+D_1L_d(q_k,q_{k+1})=0,
\end{equation}
for $k=1,\dots,N-1$. The operator $D_i$ stands for partial derivation with respect to the $i$-th argument. Later on, the notation $D_{ij} = D_j\circ D_i$ will also be used.
These equations correspond to finding critical points of the discrete action $\sum_{k=0}^{N-1}L_d(q_k,q_{k+1})$ with fixed endpoints $q_0$ and $q_N$. When discretizing a continuous Lagrangian system with Lagrangian $L \colon TQ \to \mathbb{R}$, $L_d(q_0,q_1)$ is chosen as an approximation of the action of $L$ for the solution curve joining $q_0$ to $q_1$ over a given time step $h>0$. Here, $TQ$ denotes the tangent bundle of the manifold $Q$ (the space of configurations and velocities). Refer to \cite{AbMarsdRat} for further information on tangent bundles. For further details on discrete Lagrangians and the discretization process, see for instance \cite{Marsden_West:Discrete_mechanics_and_variational_integrators} and references therein.

In \cite{MR3562389} we propose an approach to the discretization of second-order Lagrangian systems, useful for the discretization of fully actuated optimal control problems and interpolation problems. This approach consists in defining a discrete Lagrangian $L_d\colon TQ \times TQ \to \mathbb{R}$ and writing the corresponding DEL equations, with $TQ$ in the role of $Q$. As before, $L_d$ can be defined from a continuous second-order Lagrangian $L \colon T^{(2)}Q \to \mathbb{R}$, where $T^{(2)}Q$ is the second-order tangent bundle (the space of configurations, velocities and accelerations), as an approximation of the action along a solution for a given $h$.

The discrete Euler--Lagrange equations for optimal control (\textbf{DELoc}, for short) are
\begin{equation}\label{eq:DELoc}
\begin{aligned}
&D_3L_d(q_{k-1},v_{k-1},q_k,v_k)+D_1L_d(q_k,v_k,q_{k+1},v_{k+1})=0,\\
&D_4L_d(q_{k-1},v_{k-1},q_k,v_k)+D_2L_d(q_k,v_k,q_{k+1},v_{k+1})=0,
\end{aligned}
\end{equation}
for $k=1,\dots,N-1$.

Under suitable regularity conditions, these equations can be used to find a trajectory sequentially. For instance, for the DEL equations one attempts to compute $q_{k+1}$ using the previous points $q_{k-1}$ and $q_k$.
When solving boundary value problems with given initial and final conditions, some strategy should be adopted in order to arrive to the final desired condition. One such strategy is to apply a shooting method. For example, if $q_0,q_N\in Q$ and $N$ are given, one can try assigning some value to $q_1$, run the sequential algorithm and compare the resulting $q_N$ with the final condition; then adjust the value of $q_1$ and repeat the process, until the final condition is met within a certain tolerance. However, this approach can often fail to converge in practice, because of a high sensitivity of $q_N$ with respect to the starting guess $q_1$. This issue is even more prevalent in optimal control problems.

\subsection{Parallel computing}
In this paper we propose a relaxation strategy for solving boundary value problems for the discrete equations given above (see \cite{our-paper} for more details). The algorithm can be implemented using a parallel computing approach, which can significantly improve its performance and simplify the way to find approximate solutions.

We explore a strategy for simulations involving parallel computing on multicore CPUs and GPUs (graphics cards). The cores in graphics cards are processing units that are simpler and slower than regular CPU cores. However, their number presently ranges from hundreds to thousands of cores per card. This allows for great performance gains via parallelization. The approach developed here is scalable in the sense that once the algorithm for a given problem is written and tested on a GPU, additional or more powerful cards can be used to increase the number of cores and improve performance without changing the code.

The methods can be applied to problems in robotics and optimal control, by incorporating real-time feedback and trajectory correction, accounting for external perturbations or changes in the final endpoint conditions. In addition, interpolation problems can also be treated using this parallel approach.
In particular, in this paper, we will pay special attention to several types of navigation problems with prescribed boundary value conditions. 

\section{Parallel approach to the solution of the discrete equations}\label{sec:parallel_approach}
Consider the DEL equations \eqref{eq:DEL}; throughout the paper we assume all discrete Lagrangians to be at least $C^2$. Given $N\in \mathbb{N}$, $N\geq 2$, and given $q_0,q_N \in Q$, we want to find a sequence $\{q_k^*\}\equiv\{q_k^*\}_{k=0}^{N}$, with $q_0^*=q_0$, $q_N^*=q_N$, that is a solution of \eqref{eq:DEL}. Our method starts with a sequence $\{q_k\}$ chosen as initial guess, solely required to satisfy the boundary conditions, and produces a new sequence $\{\bar q_k\}$ with $\bar q_0=q_0$ and $\bar q_N=q_N$. In general, neither $\{q_k\}$ nor $\{\bar q_k\}$ will be a solution of \eqref{eq:DEL}, but by iterating this procedure we can approach a solution $\{q_k^*\}$, under certain reasonable assumptions. Sufficient conditions for convergence will be discussed in a forthcoming paper \citep{our-paper}.

\begin{figure}
	\centering
	\includegraphics{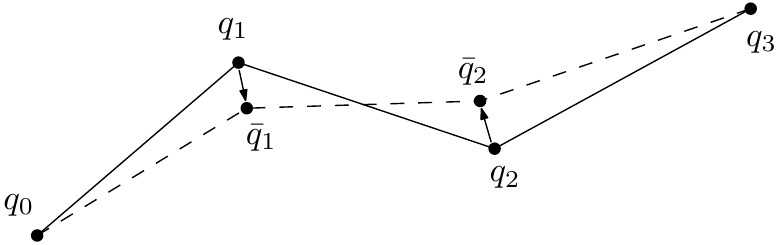}
	\caption{An iteration of the parallel method, for $N=3$.}
	\label{fig:parallel_approach}
\end{figure}

For each $k=1,\dots,N-1$, we find $\bar q_k$ by solving a modified (``parallelized'') version of \eqref{eq:DEL}:
\begin{equation}\label{eq:pDEL}
	D_2L_d(q_{k-1},\bar q_k)+D_1L_d(\bar q_k,q_{k+1})=0. 
\end{equation}
This means that for each triple $(q_{k-1}, q_k,q_{k+1})$ of points in the sequence, the middle point moves to $\bar q_k$ so that the DEL equations hold for $(q_{k-1}, \bar q_k,q_{k+1})$ (see Figure~\ref{fig:parallel_approach}). At the endpoints, we simply take $\bar q_0=q_0$ and $\bar q_N=q_N$. Computing $\bar q_k$ for all $k$ completes one iteration, and the following one will use $\{\bar q_k\}$ in place of $\{q_k\}$. This approach can be understood as the nonlinear (block) Jacobi method \citep{Vrahatis2003, Axelsson_Iterative_Solution_Methods}.

\begin{rem}
The solution $\bar q_k$ of \eqref{eq:pDEL} can be found for each $k$ independently, using the data for the neighboring points from the latest iteration. Therefore, the procedure can be performed in a parallel fashion. The computed points $\bar q_k$ remain unused until the next iteration.
\end{rem}

Assuming that $Q$ is a vector space, a less computationally expensive alternative is to replace \eqref{eq:pDEL} by a first order approximation. That is, instead of trying to solve the nonlinear system \eqref{eq:pDEL} exactly, we apply one step of the Newton-Raphson method to obtain a value for $\bar q_k$, which clearly need not coincide with the exact solution of \eqref{eq:pDEL}. This alternative update rule becomes
\begin{multline}\label{eq:NewtonRaphson}
	D_2L_d(q_{k-1}, q_{k})+D_1L_d(q_{k}, q_{k+1})
	+ \\ \left(D_{22}L_d(q_{k-1}, q_{k})+D_{11}L_d(q_{k}, q_{k+1})\right)\cdot (\bar q_k-
	q_k)=0,
\end{multline}
from which $\bar q_k$ can be computed. 
Of course, it is necessary to assume that $D_{22}L_d(q_{k-1}, q_{k})+D_{11}L_d(q_{k}, q_{k+1})$
is a regular matrix. This procedure can be seen as a single-step \textbf{Jacobi--Newton method}.

These update rules are explicit and better suited for a parallel implementation, since the same expressions can be evaluated at all time steps simultaneously, with different values for the parameters $(q_{k-1},q_k,q_{k+1})$. On the other hand, applying nonlinear solvers to \eqref{eq:pDEL} generally involves conditional statements which can cause the execution threads to diverge, that is, to execute different instructions. This can lead to a loss of performance in the parallel code.

For the DELoc equations \eqref{eq:DELoc} the procedure is analogous, since these are a particular case of the discrete Euler--Lagrange equations. Given boundary values $(q_0,v_0)$ and $(q_N,v_N)$ in $TQ$, start with an arbitrary sequence $\{(q_k,v_k)\}$ satisfying the boundary conditions. For each $k=1,\dots,N-1$, find $(\bar q_k,\bar v_k)$ by solving the parallelized equations
\begin{align*}
D_3L_d(q_{k-1},v_{k-1},\bar q_k,\bar v_k)+D_1L_d(\bar q_k,\bar v_k,q_{k+1},v_{k+1})&=0,\\
D_4L_d(q_{k-1},v_{k-1},\bar q_k,\bar v_k)+D_2L_d(\bar q_k,\bar v_k,q_{k+1},v_{k+1})&=0,
\end{align*}
and define $(\bar q_0,\bar v_0)=(q_0,v_0)$ and $(\bar q_N,\bar v_N)=(q_N,v_N)$ to get the full updated sequence $\{(\bar q_k,\bar v_k)\}$.
Iterate this procedure to approach a solution to the optimal control problem. The explicit Newton-Raphson update rule can be formulated analogously as
\begin{multline*}
\begin{bmatrix}
0 \\ 0
\end{bmatrix}= \begin{bmatrix}
D_3L_d|_{k-1}+D_1L_d|_{k}\\D_4L_d|_{k-1}+D_2L_d|_{k}
\end{bmatrix}
+\\
\begin{bmatrix}
D_{33}L_d|_{k-1}+D_{11}L_d|_{k} & D_{34}L_d|_{k-1}+D_{12}L_d|_{k} \\D_{43}L_d|_{k-1}+D_{21}L_d|_{k} & D_{44}L_d|_{k-1}+D_{22}L_d|_{k}
\end{bmatrix}
\begin{bmatrix}
\bar q_k-q_k \\ \bar v_k-v_k
\end{bmatrix},
\end{multline*}
where we have used the shorthand notation $D_1L_d|_k\equiv D_1L_d(q_{k},v_{k}, q_{k+1}, v_{k+1})$, and similarly for the other derivatives. As before, we need to assume regularity of the coefficient matrix.

\section{Zermelo's navigation problem}

Zermelo's navigation problem \citep{Zermelo, Bao, Kopacz, java} is usually presented as a time-optimal control problem, which aims to find the minimum time trajectories on a Riemannian manifold $(Q, g)$ under the influence of a drift vector field $W\in {\mathfrak X}(Q)$, which can be interpreted as wind (or water currents). The goal is to navigate from a point to another in $Q$ along a path $\gamma(s)$ in the shortest possible time in the presence of this wind, assuming that the ship engine provides a constant thrust relative to it, that is, $|\dot{\gamma}(s)-W(\gamma(s))|=1$, where $|\cdot|$ denotes the norm provided by $g$. It is assumed that $|W(q)|< 1$ for all $q\in Q$. 

These minimum time trajectories are precisely the geo\-de\-sics for a particular type of Finsler metric, a Randers metric defined by (see \cite{Bao} and references therein)
\begin{equation*}
F(q,v_q)=\sqrt{a(v_q, v_q)}+\langle b(q), v_q\rangle
\end{equation*}
where 
\begin{align*}
	a(v_q, v_q)&=\frac{1}{\alpha(q)}g(v_q,v_q)+\frac{1}{\alpha(q)^2}g(W(q), v_q)^2\\
  \langle b(q), v_q\rangle&=-\frac{1}{\alpha(q)} g(W(q), v_q)=-\left\langle \frac{\flat_g({W}(q))}{\alpha(q)}, v_q\right\rangle \\
  \alpha(q)&=1-g(W(q),W(q))=1-|W(q)|^2>0.
\end{align*}
Here $\flat_g\colon {\mathfrak X}(Q)\rightarrow \Omega^1(Q)$ is the musical isomorphism defined by $\langle \flat_g(X), Y\rangle=g(X,Y)$ for all $X, Y\in {\mathfrak X}(Q)$. 

The time it takes the ship to move along a curve $\gamma\colon [t_0, t_N]\to Q$ is
\begin{equation}\label{eq:intF}
  \int_{t_0}^{t_N} F(\gamma(s), \dot \gamma(s))\,ds.
\end{equation}
Note that this integral is invariant under orientation-preserving reparametrizations of $\gamma$, since Finsler metrics are positively 1-homogeneous, that is, $F(q,\lambda v_q)=\lambda F(q,v_q)$ for any $\lambda>0$. Therefore, the solution curves are not unique. In fact, $F$ is not regular as a Lagrangian function. Similar to the case of Riemannian metrics and the problem of minimizing length or energy, this can be circumvented by considering instead the functional
\begin{equation}\label{eq:intF2}
  \int_{t_0}^{t_N} \big(F(\gamma(s),\dot\gamma(s))\big)^2\,ds .
\end{equation}
Any extremal of this functional will be an extremal of \eqref{eq:intF}, and any extremal of \eqref{eq:intF} admits an orientation-preserving reparametrization that makes it an extremal of \eqref{eq:intF2} (see \cite{Masiello} and references therein).

As a particular case, consider $Q={\mathbb R}^2$ with the Euclidean metric, where we are to find critical curves $(x,y)=(x(s),y(s))$ for the functional
\begin{align*}
\int_{t_0}^{t_N}&\left[\sqrt{\frac{1}{\alpha}(\dot{x}^2+\dot{y}^2)+\frac{1}{\alpha^2}(W_1(x,y)\dot{x}+W_2(x, y)\dot{y})^2}\right.\\
&\left.-\frac{1}{\alpha} \left(W_1(x, y)\dot{x}+W_2(x, y)\dot{y}\right)\right]^2\, ds
\end{align*}
with $\alpha=1-(W_1^2+W_2^2)$.

\begin{figure}
  \begin{center}
  \includegraphics[scale=.7]{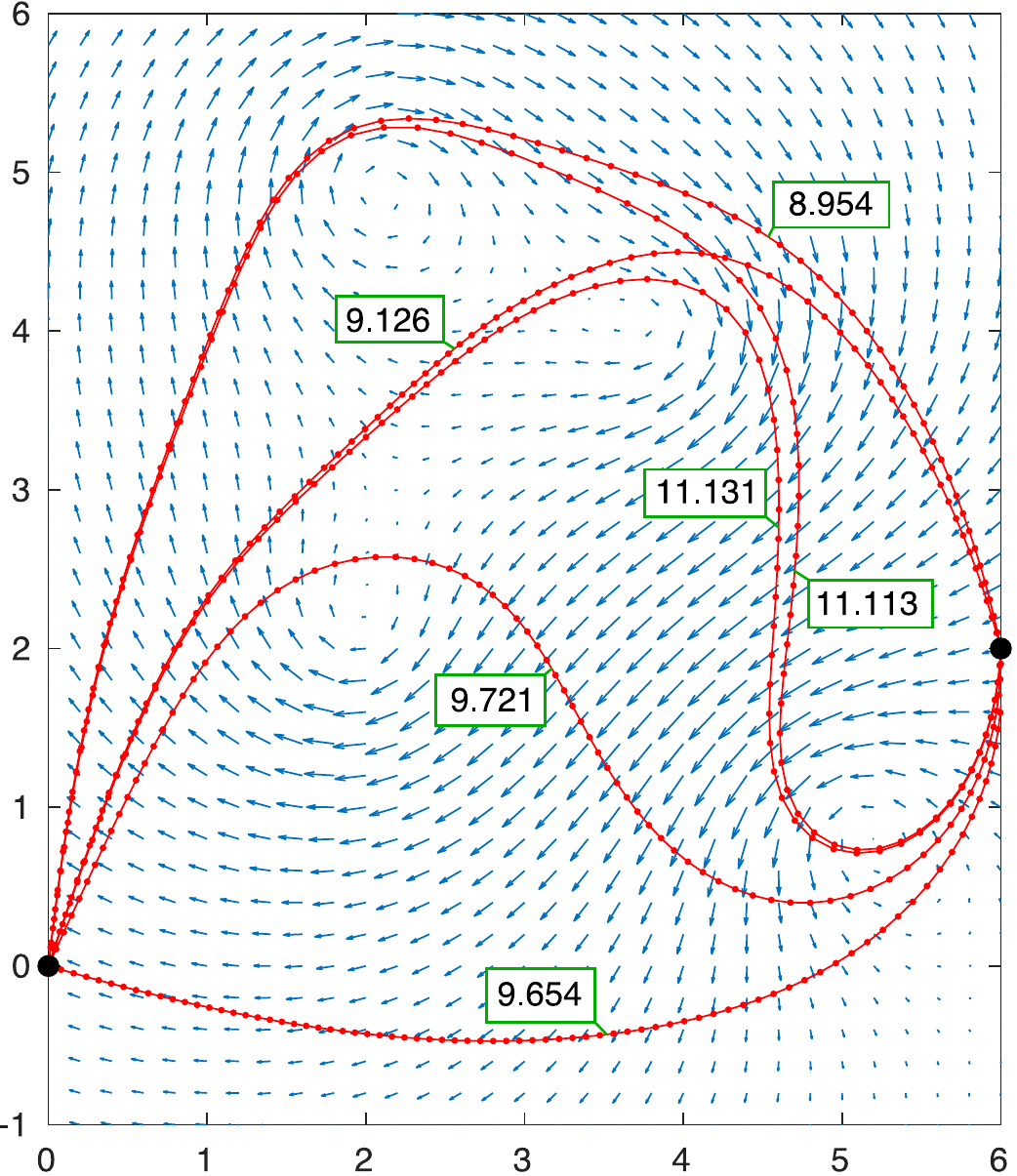}
  \end{center}
  \caption{Several local solutions to the optimal time navigation problem starting from $(0,0)$ and ending at $(6,2)$. The time for each trajectory is shown.}\label{fig:timeoptimal}
\end{figure}

Figure \ref{fig:timeoptimal} shows six local solutions to the navigation problem found using our approach, starting at $(0,0)$ and ending at $(6,2)$, for the vector field
$W=1.7\cdot(-R_{2,2}-R_{4,4}-R_{2,5}+R_{5,1})$, where
\[
  R_{a,b}(x,y)=\frac{1}{3((x-a)^2+(y-b)^2)+1}\begin{bmatrix}
    -(y-b)\\x-a
  \end{bmatrix}.
\]
 The scale factor $1.7$ was chosen so that the maximum value of $|W|$ is almost 1. These solutions were obtained using different piecewise straight lines as initial guesses for the trajectories, having between two and five line segments. In this case, $N=80$. The total navigation time~\eqref{eq:intF}, which is displayed beside each trajectory, is locally optimal. In general, finding a global minimum will require exploring different initial guesses.

 As the discrete Lagrangian, we used
 \[
   L_d(q_0,q_1)=\frac{h}{2} \left(F^2\left(q_0,\frac{q_1-q_0}{h}\right)+F^2\left(q_1,\frac{q_1-q_0}{h}\right)\right).
 \]

\section{Fuel-optimal navigation problem}\label{fuel}
We also consider a non-equivalent variant of Zermelo's problem. If $T>0$ is a fixed time, we seek trajectories minimizing the cost function
\begin{equation*}
\int_0^T \frac{1}{2}(u_1^2+u_2^2)\; dt\,,
\end{equation*}
which can be interpreted as a measure of fuel expenditure.
The system is subject to the control equations
\begin{align*}
	\dot{x}&= u_1 + W_1(x, y)\,,\\
	\dot{y}&= u_2 + W_2(x, y)\,.
\end{align*}
The goal is to arrive at a given destination at time $T$, extremizing fuel expenditure with no a priori bounds on the engine's power. 
This problem is equivalent to solving the Euler-Lagrange equations for the Lagrangian
\begin{equation*}
L(x, y, \dot{x}, \dot{y}) = \frac{1}{2} \left[(\dot{x}-W_1(x,y))^2+ (\dot{y}-W_2(x,y))^2\right]\,,
\end{equation*}
with fixed $(x(0), y(0))$ and $(x(T), y(T))$ as boundary conditions.

For our simulations, we considered $W(x,y)=(
  \cos(2x-y-6),
  \frac{2}{3}\sin(y)+x-3)$.
We discretized the Lagrangian as 
\[
   L_d(q_0,q_1)=\frac{h}{2} \left(L\left(q_0,\frac{q_1-q_0}{h}\right)+L\left(q_1,\frac{q_1-q_0}{h}\right)\right).
 \]

\begin{figure}
  \begin{center}
  \includegraphics[scale=.7]{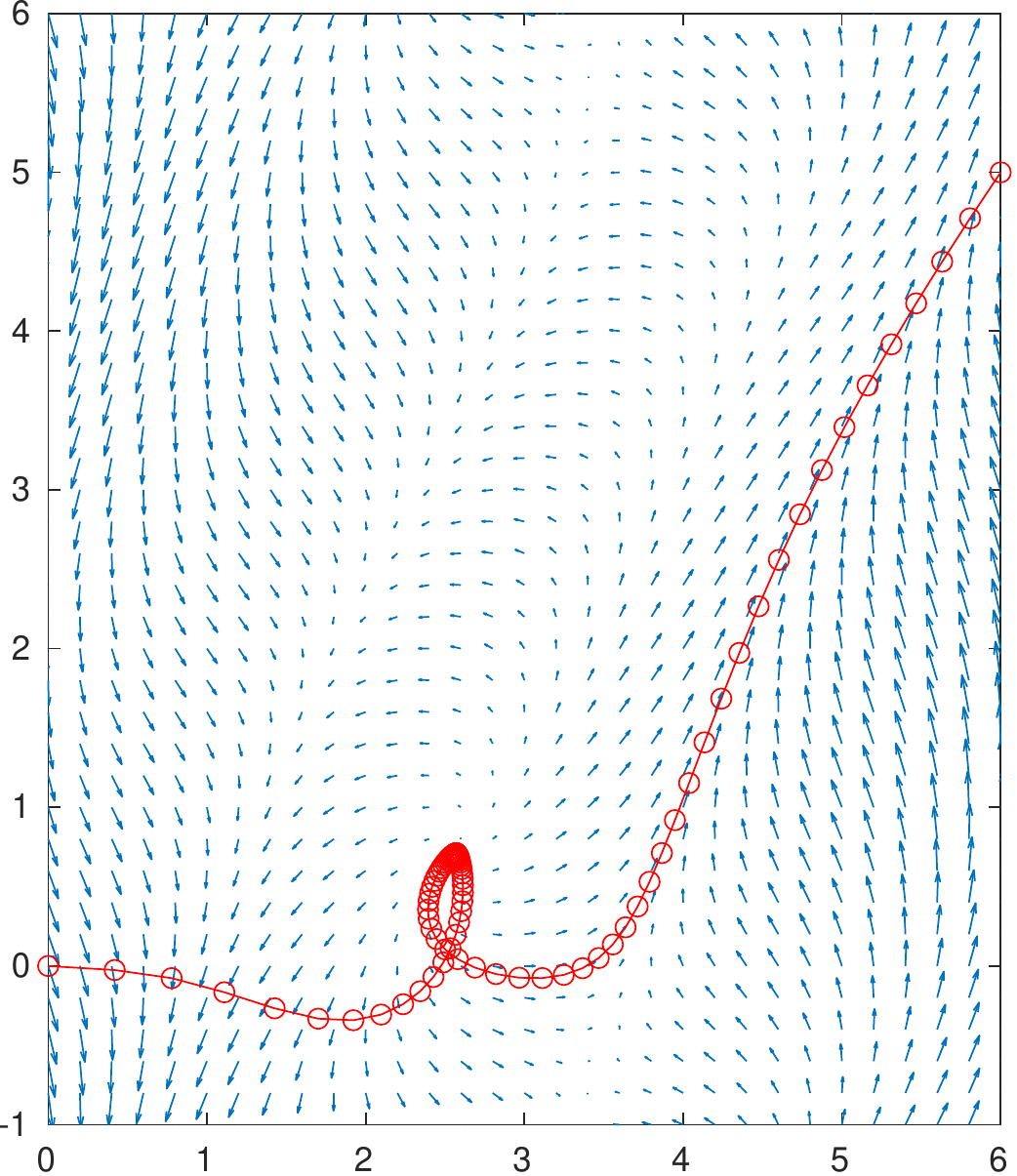}
  \end{center}
  \caption{A minimal fuel trajectory for a fixed total duration $T=30$, joining $(0,0)$ to $(6,5)$, with $N=200$. }\label{fig:fueloptimal}
\end{figure}

Figure \ref{fig:fueloptimal} shows a trajectory that has a locally optimal fuel expenditure among the discrete curves joining the given points $(0,0)$ and $(6,5)$ in time $T=30$. A straight line was used as the initial guess. Notice that in the first part of its journey, the ship travels to an equilibrium point of $W$, where it awaits the right moment to continue to its destination, which must be reached at the specified time. We emphasize that in this variant of the problem the total travel time is imposed externally. Other values for $T$ will have optimal trajectories with different fuel expenditure. For $T<9$ (approximately) the optimal trajectories do not pass near the equilibrium mentioned above.

\section{Interpolation problems}
In this section we explore the application of our parallel iterative method to the case of a second-order Lagrangian system in the context of interpolation problems. In a biomedical setting, this kind of problem appears when comparing a series of images in longitudinal studies \citep{Crouch,Invariant1,Invariant2}. Observe that the extension to higher-order Lagrangians is straightforward. 
First, consider a second-order Lagrangian $L\colon T^{(2)}Q\to {\mathbb R}$ and let $N\in \mathbb{N}$, $[t_0,t_N]\subset \mathbb{R}$.
We say that a curve $q\colon [t_0, t_N]\to Q$ is {\bf critical} for the action
\begin{equation} \label{ElasticSplines}
\mathcal{J} [q]:= \int_{t_0}^{t_N} L(q(t), \dot{q}(t), \ddot{q}(t))\; dt
\end{equation}
among all curves with given values and first derivatives fixed at the endpoints, i.e., 
\[q(t_0), \dot{q}(t_0), q(t_N), \dot{q}(t_N)\qquad \hbox{(boundary conditions)}
\]
if and only if $q$ is a solution of the second order Euler-Lagrange equations given by the system of fourth-order differential equations: 
\begin{equation}\label{second}
\frac{d^2}{dt^2}\left(
\frac{\partial L}{\partial \ddot{q}}\right)-\frac{d}{dt}\left(
\frac{\partial L}{\partial \dot{q}}\right)+\frac{\partial L}{\partial q}=0
\end{equation}
Additionally, we assume that we have $l+1$ interpolation points or knots $\hat{q}_a\in Q$, $a=0,\dots ,l$, which are reached at times 
$\hat{t}_a =t_0+ N_a h$,
where $h=(t_N-t_0)/N$, $N_a\in \{0, 1, \ldots, N\}$,
$N_a< N_b$ if $a<b$, with $N_0=0$ and $N_l=N$.

In order to discretize this problem, we replace the integral~\eqref{ElasticSplines} by a sum over times $t_k=t_0+kh$ for $k=0, \ldots, N$. 
Following our approach in \cite{MR3562389}, we discretize the action as
\begin{equation} \label{ElasticSplines-d}
\mathcal{J}_d := \sum_{k=0}^{N-1} L_d(q_k, v_k, q_{k+1}, v_{k+1})
\end{equation}
where $L_d\colon TQ\times TQ\rightarrow {\mathbb R}$ is a discretization of $L$.
Moreover, the problem is subject to the interpolation constraints
\begin{equation}\label{aq1}
q_{N_a}=\hat{q}_a , \quad\text{for all }a=1,\ldots, l-1
\end{equation}
and the boundary conditions
\begin{equation}\label{aq2}
q_0=\hat{q}_0,\quad v_0= \hat{v}_0 ,\quad\text{and}\quad
q_N=\hat{q}_l,\quad v_N= \hat{v}_l .
\end{equation}

Our parallel integrator works as follows. 
Take an arbitrary sequence $\{(q_k, v_k)\}$ satisfying the interpolation constraints (\ref{aq1}) and the boundary conditions (\ref{aq2}).
Now construct the sequence $\{(\bar{q}_k, \bar{v}_k)\}$ by solving the parallelized problem 
\begin{align*}
D_3L_d(q_{k-1}, v_{k-1}, \bar{q}_{k}, \bar{v}_k)
+D_1L_d(\bar{q}_{k}, \bar{v}_k ,q_{k+1}, v_{k+1})&=0\; ,\\
D_4L_d(q_{k-1}, v_{k-1}, \bar{q}_{k}, \bar{v}_k)
+D_2L_d(\bar{q}_{k}, \bar{v}_k ,q_{k+1}, v_{k+1})&=0
\end{align*}
if $1\leq k\leq {N}-1$ and $k\not=N_a$, $1\leq a \leq l - 1$. At each knot $k=N_a$, $1\leq a\leq l-1$, take $\bar{q}_{N_a}=q_{N_a}$ and compute $\bar{v}_{N_a}$ by solving the equation
\begin{multline*}
D_4L_d(q_{N_a-1}, v_{N_a-1}, \bar{q}_{N_a}, \bar{v}_{N_a})\\
+D_2L_d(\bar{q}_{N_a}, \bar{v}_{N_a}, q_{N_a+1}, v_{N_a+1})=0\; .
\end{multline*}
Finally, take $(\bar{q}_0,\bar{v}_0)=(q_0,v_0)$, $(\bar{q}_N,\bar{v}_N)=(q_N,v_N)$. 
Observe that the derived sequence $\{(\bar{q}_k, \bar{v}_k)\}$, $k=0, \ldots, N$
satisfies the 
interpolation constraints
\[
\bar{q}_{N_a}=\hat{q}_a , \quad\text{for all } a=1,\ldots, l-1
\]
and the boundary conditions 
\begin{equation*}
\bar q_0=\hat{q}_0,\quad \bar v_0= \hat{v}_0 ,\quad\text{and}\quad
\bar q_N=\hat{q}_l,\quad \bar v_N= \hat{v}_l \; .
\end{equation*}
By iterating this procedure, we approach a trajectory having a locally optimum value of the cost functional (\ref{ElasticSplines-d}).

\subsection{An application: fuel-optimal control problem with a weight minimizing the total variation in the control variables }
As a modification of the application given in Section \ref{fuel}, consider the following optimal control problem. Our aim is still minimizing the fuel expenditure functional while also minimizing the total variation in the control variables. 
Now the goal is to navigate from a departure point to a destination point passing through given waypoints (knots) at prescribed times, minimizing the cost functional
\[
\int_0^T \frac{1}{2}(u_1^2+u_2^2 +c v_1^2+cv_2^2)\; dt
\]
subject to the control equations
\begin{equation*}
\begin{array}{rclrcl}
\dot{x} &=& u_1+W_1(x, y), & \quad\dot{y} &=& u_2+W_2(x,y),\\
\dot{u}_1&=&v_1, & \dot{u}_2&=&v_2\,.
\end{array}
\end{equation*}
Here $c > 0$ is a weight.

The continuous problem is equivalent to solving the fourth-order Euler-Lagrange equations for the second-order Lagrangian
\begin{multline*}
	L(x, y, \dot{x}, \dot{y}, \ddot{x}, \ddot{y})=\frac{1}{2} \left[(\dot{x}-W_1(x,y))^2+ (\dot{y}-W_2(x,y))^2\right.\\	
	+c\left(\ddot{x}- D_1W_1(x(t),y(t))\dot{x}- 
	D_2W_1(x(t),y(t))\dot{y}\right)^2
	\\
	\left.
	+c\left(\ddot{y}- D_1W_2(x(t),y(t))\dot{x}- 
	D_2W_2(x(t),y(t))\dot{y}\right)^2
	\right]
\end{multline*}
As boundary conditions, we consider $(q(0), \dot q(0))$ and $(q(T), \dot q(T))$ fixed. In addition, the system is subject to the
interpolation constraints
\begin{equation}\label{aq1example}
q(\hat{t}_a)=\hat{q}_a , \quad\text{for all }a=1,\ldots, l-1
\end{equation}
with $0<\hat{t}_a<\hat{t}_b<T$ for all $a, b\in\{1,\ldots, l-1\}$ and $a<b$.

As a discretization of the cost function we propose, for instance, a 2-stage Lobatto discretization:
\begin{align*}
&L_d(q_k, v_k, q_{k+1}, v_{k+1}) =\\
&\quad\frac{h}{2} \left[ L\left(q_k, v_k, \frac{2}{h^2}(3(q_{k+1}-q_k)-h(v_{k+1}+2v_k)\right) \right.\\
&\quad + \left. L\left(q_{k+1}, v_{k+1}, -\frac{2}{h^2}(3(q_{k+1}-q_k)-h(2v_{k+1}+v_k)\right)\right]
\end{align*}

\begin{figure}
  \begin{center}
  \includegraphics[scale=.7]{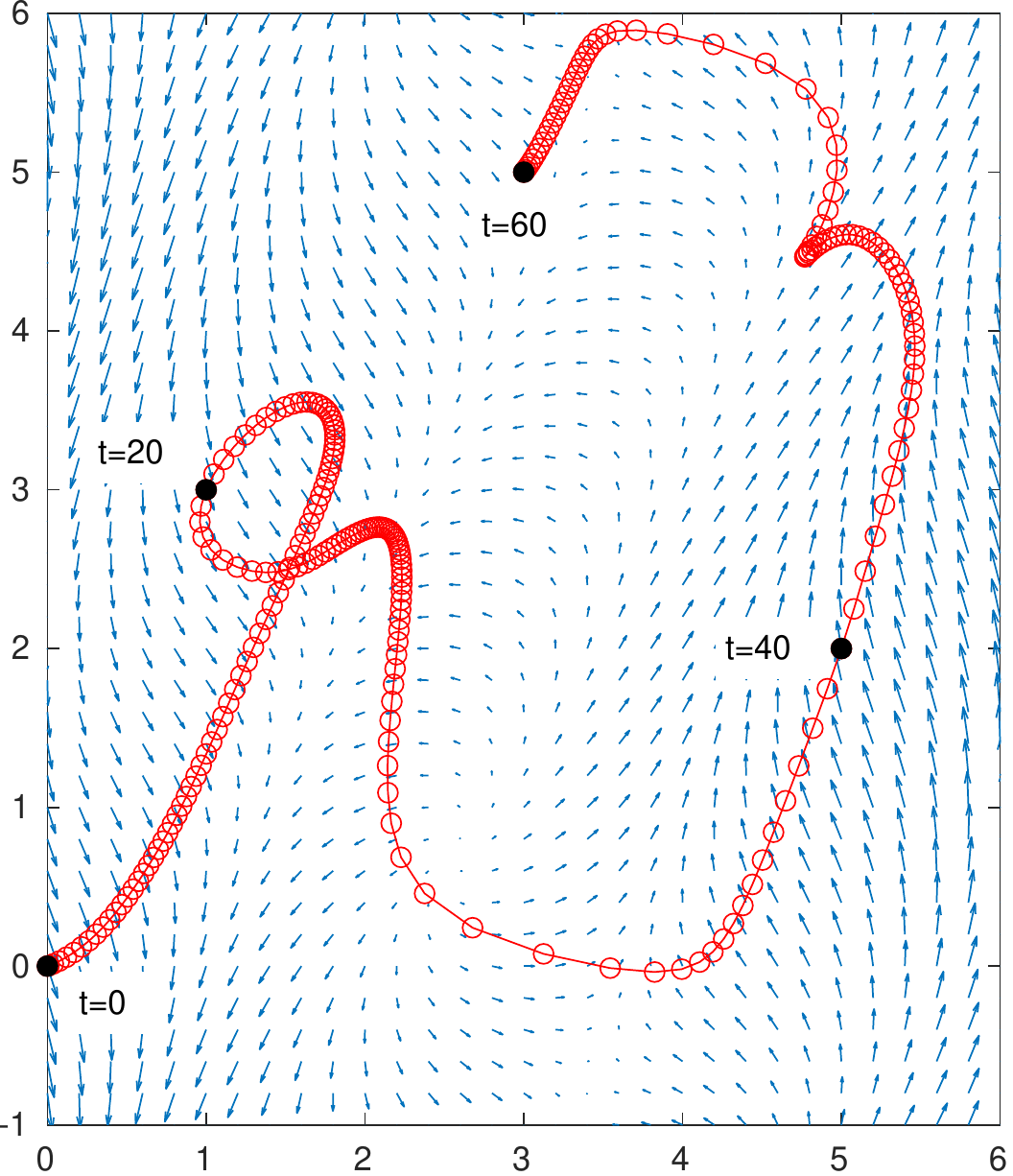}
  \end{center}
  \caption{An optimal trajectory for the second-order problem with interpolation nodes.}\label{fig:secondorder}
\end{figure}

In Figure \ref{fig:secondorder} we show an optimal trajectory starting at $(0,0)$ and ending at $(3,5)$ at $T=60$, with zero velocity at both endpoints, and passing through the prescribed positions $(1,3)$ and $(5,2)$ at times $20$ and $40$ respectively. The vector field $W$ is the same as in the previous example. We used $c=50$ and $N=240$. The initial guess was a cubic spline connecting the four waypoints and having the prescribed velocity at the endpoints. If a shooting strategy was applied, a small change in $(q_1,v_1)$ would produce a very different trajectory.

All of the examples have been computed using the Jacobi-Newton method mentioned in section~\ref{sec:parallel_approach}. In this last simulation, a ``damping'' coefficient $\delta=.05$ was used as follows. If $(\Delta q_k,\Delta v_k)$ is the increment computed by the Jacobi-Newton method at the index $k$, then the update rule applied was $(\bar{q}_k,\bar{v}_k)=(q_k,v_k)+(1-\delta)(\Delta q_k,\Delta v_k)$. This was used to control an instability of the trajectory that appears around $t=50$.

\section{Comments on the implementation}
Prototypes of the algorithms were implemented using MATLAB\textsuperscript{\circledR}
 and the Parallel Computing Toolbox\textsuperscript{TM} (\cite{matlabparallel}). The simulations were run on an Intel\textsuperscript{\circledR} Core\textsuperscript{TM} i5-7400 CPU and an NVIDIA GeForce GT 740 GPU. The execution times provided throughout this section are given only for illustration purposes and transparency. However, these are strongly dependent on the implementation, which may have been suboptimal.
 
 For Figures 2 and 3, the stopping criterion was
\[
  \max_k{\|D_2L_d(q_{k-1}, q_k)+D_1L_d(q_k,q_{k+1})\|}<10^{-4}h^2.
\]
The factor $h^2$ is included because the left-hand side of this inequality is $O(h)$ for any smooth curve sampled at times $kh$, regardless of it being a trajectory of the system, and for $L_d$ approximating the exact discrete Lagrangian. Each one of the six trajectories in Figure 2 took between 13 to 16 seconds to compute, and required around 11000 to 13000 iterations; however, a trajectory with a time expenditure equal (within three decimal places) to the optimal time shown in the figure was already obtained in the first 1 or 2 seconds of each simulation. Setting the tolerance to $10^{-10}h^2$ does not decrease the time expenditure within three decimal places.

For Figure 3, after 227000 iterations (about 237 seconds) the trajectory fulfilled the stopping criterion, with a fuel expenditure of $5.597$. As before, a solution with the same fuel expenditure within three decimal places one was already obtained after 22000 iterations (23 seconds). This relatively long time was due to the fact that the initial guess, regularly spaced points along a  straight line segment, was not close to the solution.

Finally, the interpolation problem (Figure 4) was significantly slower to converge. The initial guess was a cubic spline sampled with $N=120$. Although a natural choice, a spline is not a high-quality guess as the actual solution is a quite convoluted curve. After 150000 iterations (around 480 seconds), the trajectory started to look qualitatively like the final solution, with a total fuel expenditure of $134.2$. After this occurred, we refined the trajectory by taking $N=240$ and halving the time step. Letting the simulation run for 5000 more seconds, we obtained a slightly better trajectory, with fuel expenditure $133.4$.

\section{Conclusions and future work}

In this paper we have shown that discrete variational methods combined with a parallel iterative approach are well-suited for boundary value problems, where techniques such as the shooting method would fail. This has been tested in three examples related with navigation problems: the classical minimum-time Zermelo's navigation problem, a fuel-optimal navigation problem consisting in arriving to the destination at a fixed time and, finally, an interpolation problem for a fuel-optimal control problem with a weight minimizing the total variation in the control variables. The numerical simulations give us trajectories which locally minimize the corresponding cost (time, fuel expenditure, etc). 

In \cite{our-paper} we prove rigorous conditions for the convergence of Equations (\ref{eq:pDEL}) and (\ref{eq:NewtonRaphson}) and also of the corresponding extensions to second-order Lagrangians. In that paper we will also show that time-dependent water currents  can be added to the proposed  navigation problems. In a forthcoming paper  we will describe how to adapt our parallel iterative method for the case of invariant first and second order Lagrangian systems where the configuration space is a Lie group. 

Problems involving sharp univariate constraints such as state space or control space exclusion zones have not yet been studied using this approach. Nevertheless, we believe that such problems may be handled by means of penalty potentials. Moreover, the iterative procedure of increasing the penalty to approximate sharp boundaries could be coupled with the iteration of our algorithm. We intend to address this in a forthcoming paper.

\bibliography{references}             

\newcommand\oneletter[1]{#1}\newcommand\Yu{Yu}
\begin{thebibliography}{15}
\providecommand{\natexlab}[1]{#1}
\providecommand{\url}[1]{\texttt{#1}}
\providecommand{\urlprefix}{URL }
\expandafter\ifx\csname urlstyle\endcsname\relax
  \providecommand{\doi}[1]{doi:\discretionary{}{}{}#1}\else
  \providecommand{\doi}{doi:\discretionary{}{}{}\begingroup
  \urlstyle{rm}\Url}\fi

\bibitem[{Abraham et~al.(1988)Abraham, Marsden, and Ratiu}]{AbMarsdRat}
Abraham, R., Marsden, J.E., and Ratiu, T. (1988).
\newblock \emph{Manifolds, tensor analysis, and applications}, volume~75 of
  \emph{Applied Mathematical Sciences}.
\newblock Springer-Verlag, New York, second edition.

\bibitem[{Axelsson(1994)}]{Axelsson_Iterative_Solution_Methods}
Axelsson, O. (1994).
\newblock \emph{Iterative solution methods}.
\newblock Cambridge University Press, Cambridge.

\bibitem[{Bao et~al.(2004)Bao, Robles, and Shen}]{Bao}
Bao, D., Robles, C., and Shen, Z. (2004).
\newblock Zermelo navigation on {R}iemannian manifolds.
\newblock \emph{J. Differential Geom.}, 66(3), 377--435.

\bibitem[{Colombo et~al.(2016)Colombo, Ferraro, and Mart\'\i n~de
  Diego}]{MR3562389}
Colombo, L., Ferraro, S., and Mart\'\i n~de Diego, D. (2016).
\newblock Geometric integrators for higher-order variational systems and their
  application to optimal control.
\newblock \emph{J. Nonlinear Sci.}, 26(6), 1615--1650.

\bibitem[{Crouch and Silva~Leite(1995)}]{Crouch}
Crouch, P. and Silva~Leite, F. (1995).
\newblock The dynamic interpolation problem: on {R}iemannian manifolds, {L}ie
  groups, and symmetric spaces.
\newblock \emph{J. Dynam. Control Systems}, 1(2), 177--202.

\bibitem[{Ferraro et~al.(2021)Ferraro, Mart{\'\i}n~de Diego, and Sato
  Mart\'{\i}n~de Almagro}]{our-paper}
Ferraro, S., Mart{\'\i}n~de Diego, D., and Sato Mart\'{\i}n~de Almagro, R.T.
  (2021).
\newblock A parallel iterative method for variational integration.
\newblock \emph{Work in progress}.

\bibitem[{Gay-Balmaz et~al.(2012{\natexlab{a}})Gay-Balmaz, Holm, Meier, Ratiu,
  and Vialard}]{Invariant1}
Gay-Balmaz, F., Holm, D.D., Meier, D.M., Ratiu, T.S., and Vialard, F.X.
  (2012{\natexlab{a}}).
\newblock Invariant higher-order variational problems.
\newblock \emph{Comm. Math. Phys.}, 309(2), 413--458.

\bibitem[{Gay-Balmaz et~al.(2012{\natexlab{b}})Gay-Balmaz, Holm, Meier, Ratiu,
  and Vialard}]{Invariant2}
Gay-Balmaz, F., Holm, D.D., Meier, D.M., Ratiu, T.S., and Vialard, F.X.
  (2012{\natexlab{b}}).
\newblock Invariant higher-order variational problems {II}.
\newblock \emph{J. Nonlinear Sci.}, 22(4), 553--597.

\bibitem[{Javaloyes and S\'{a}nchez(2017)}]{java}
Javaloyes, M.A. and S\'{a}nchez, M. (2017).
\newblock Wind {R}iemannian spaceforms and {R}anders-{K}ropina metrics of
  constant flag curvature.
\newblock \emph{Eur. J. Math.}, 3(4), 1225--1244.

\bibitem[{Kopacz(2019)}]{Kopacz}
Kopacz, P. (2019).
\newblock On generalization of {Z}ermelo navigation problem on {R}iemannian
  manifolds.
\newblock \emph{Int. J. Geom. Methods Mod. Phys.}, 16(4), 1950058, 19.

\bibitem[{Marsden and
  West(2001)}]{Marsden_West:Discrete_mechanics_and_variational_integrators}
Marsden, J.E. and West, M. (2001).
\newblock Discrete mechanics and variational integrators.
\newblock \emph{Acta Numer.}, 10, 357--514.

\bibitem[{Masiello(2009)}]{Masiello}
Masiello, A. (2009).
\newblock An alternative variational principle for geodesics of a {R}anders
  metric.
\newblock \emph{Adv. Nonlinear Stud.}, 9(4), 783--801.

\bibitem[{Math\-Works\textsuperscript{\circledR}(2019)}]{matlabparallel}
Math\-Works\textsuperscript{\circledR} (2019).
\newblock \emph{MATLAB\textsuperscript{\circledR} and Parallel Computing
  Toolbox\textsuperscript{TM}}.
\newblock Natick, Massachusetts, United States.

\bibitem[{Vrahatis et~al.(2003)Vrahatis, Magoulas, and
  Plagianakos}]{Vrahatis2003}
Vrahatis, M.N., Magoulas, G.D., and Plagianakos, V.P. (2003).
\newblock From linear to nonlinear iterative methods.
\newblock \emph{Appl. Numer. Math.}, 45(1), 59--77.
\newblock 5th IMACS Conference on Iterative Methods in Scientific Computing
  (Heraklion, 2001).

\bibitem[{Zermelo(1931)}]{Zermelo}
Zermelo, E. (1931).
\newblock Über das {N}avigationsproblem bei ruhender oder veränderlicher
  {W}indverteilung.
\newblock \emph{Z. Angew. Math. Mech.}, 11, 114--124.

\end{thebibliography}
                                                   
\end{document}